\documentclass[11pt]{article}

\usepackage[dvips]{graphics}

\usepackage{amsmath,amsfonts,amssymb,amsthm}
\usepackage{xspace}
\usepackage{epsfig}
\usepackage{psfrag}
\usepackage{equation}
\usepackage{a4wide}
\newcommand{\ip}[1]{\langle #1 \rangle}

\newcommand{\mat}{\textsc{Matlab}\xspace}

\newcommand{\lemref}[1]{{Lemma \ref{#1}}}

\newtheorem{theorem}{Proposition}
\newtheorem{lemma}{Lemma}

\def\Z{\mathbb Z}
\def\R{\mathbb R}

\begin{document}

\title{Finite to infinite steady state solutions, bifurcations of an integro-differential equation}

\author{S. K. Bhowmik$^1$, D. B. Duncan$^2$, M. Grinfeld$^3$, G. J. Lord$^4$}

\maketitle


{\footnotesize
 \centerline{$^1$KdV Institute for Mathematics, University of Amsterdam,
   Amsterdam, NL, S.K.Bhowmik@uva.nl.}
 \centerline{$^2$Department of Mathematics and Maxwell Institute,
   Heriot--Watt University,Edinburgh, UK, D.B.Duncan@hw.ac.uk}
 \centerline{$^3$Department of Mathematics, University of Strathclyde,
   Glasgow, UK, M.Grinfeld@strath.ac.uk}
 \centerline{$^4$Department of Mathematics and Maxwell Institute,
   Heriot--Watt University,Edinburgh, UK, G.J.Lord@hw.ac.uk}
}

\begin{abstract}
We consider a bistable integral equation which governs the stationary
solutions of a convolution model of solid--solid phase transitions on
a circle. We study the bifurcations of the set
of the stationary solutions as the diffusion coefficient is varied to
examine the transition from an infinite number of steady states to
three for the continuum limit of the semi--discretised system.   
We show how the symmetry of the problem is responsible for the generation
and stabilisation of equilibria and comment on the puzzling connection
between continuity and stability that exists in this problem. 
\end{abstract}

\section{Introduction}\label{introd}

Integro-differential equations are used to model various
phenomena in materials science~\cite{P.Bates,A.Chmaj,C.Fife,
  F.Chen,Dug,HrtlyWnnr} and biology~\cite{J.C.Louis,
  Deng,J.Medlock,WilsonCowan}, which involve non-local diffusion/dispersal
mechanisms. We consider the integro-differential equation (IDE)
\begin{equation}\label{eq:ut}
u_{t}=\varepsilon \left( \int_{\R}\ J^{\infty}(x-y) u(y,t)dy-u(x,t)
\int_{\R}\ J^{\infty}(x-y)dy\right)+f(u),
\end{equation}
where the $L^1(\R)$ kernel $J^\infty$ satisfies
$J^{\infty}(x)\ge 0,$ $J^{\infty}(x)=J^{\infty}(-x)$ and $f(u)$ is a
bistable nonlinearity. 
Below we routinely consider $f(u)=u(1-u^2)$ and kernel 
\begin{equation}\label{eq:expk}
J^\infty(x) = \sqrt{\frac{100}{\pi}} \exp(-100 x^2),
\end{equation}
so that $\int_\R J^\infty \, \ dx = 1$.
To obtain a well-defined problem, (\ref{eq:ut}) has to be supplemented by a
suitable initial condition,
$u(x,0)=u_0(x)$ which needs to be chosen in a suitable function space, see
\cite{GHHMV,vhmg01,GS06}.

The convolution equation (\ref{eq:ut}) is the $L^2$-gradient flow of
the free energy functional
\begin{equation}\label{fe:f}
E(u)=\frac{1}{4}\varepsilon
\int_{\R}\int_{\R}J^{\infty}(x-y)\left(u(y)-u(x)\right)^2dxdy+
\int_{\R}\hat{F}(u)dx,
\end{equation}
where $\hat{F}(u,t)$ is the smooth double well potential, $
\hat{F}'(u) = -f(u)$.

For an overview of the use of (\ref{eq:ut}) in materials science, see
\cite{phase03}. There are many papers dealing with the mathematical
analysis of this equation, which examine existence and stability of
travelling waves \cite{C.Fife}, the structure of the stationary
solutions set \cite{A.Chmaj}, propagation of discontinuities
\cite{CFphase}, coarsening \cite{Dug} and long time behaviour
\cite{GHHMV,vhmg01,Rossi07,Rossi09}. 

Note, in particular, that in \cite{GHHMV} it is shown that if the diffusion
coefficient $\varepsilon$ is sufficiently large, a ``Conway--Hopf--Smoller''
type result holds: the only stable steady state solutions, say, in
$L^\infty(\R)$, are the constant stable steady states of the kinetic
equation $u_t=f(u)$. Thus, if we choose $f(u)=u(1-u^2)$, the stable states are
$u=1$ and $u=-1$. On the other hand, if $\varepsilon=0$, (\ref{eq:ut})
admits an uncountable set of equilibria: let $X$, $Y$ and $Z$ be any
disjoint sets such that $A\cup B \cup C= \R$, then a function $u(x)$ that is
equal to $1$ on $X$, $-1$ on $Y$ and $0$ on $Z$ is a steady state solution.
Note that if $Z = \emptyset$, all the resulting equilibria are stable in
$L^\infty (I)$. Furthermore, it is shown in \cite{Dug} that there exists an
$\varepsilon_0>0$ which depends on the kernel $J^\infty$, such that for all
$0< \varepsilon < \varepsilon_0$ the set of steady state solutions of
(\ref{eq:ut}) is in one-to-one correspondence with the set of equilibria of
$u_t=f(u)$. Hence, in view of the above, it is of interest to perform a
bifurcation analysis of the set of steady states of
(\ref{eq:ut}),
\begin{equation}\label{eq:ss}
0 = \varepsilon \int_\R J^\infty(x-y)(u(y)-u(x))\, dy + f(u),
\end{equation}
as we decrease $\varepsilon$ from some initially large value to zero,
and investigate the transition from a finite to infinite set of solutions.

To the best of our knowledge, such a study has not been performed before.
The object of this paper is precisely such a study of the spatially
discretised version of (\ref{eq:ut}).  For simplicity, here we restrict
ourselves to $1$-periodic patterns.

\begin{figure}[!h]
\begin{center}
{\bf (a) \hspace{0.44\textwidth} (b)}
\includegraphics[width=0.89\textwidth,height=0.22\textheight]{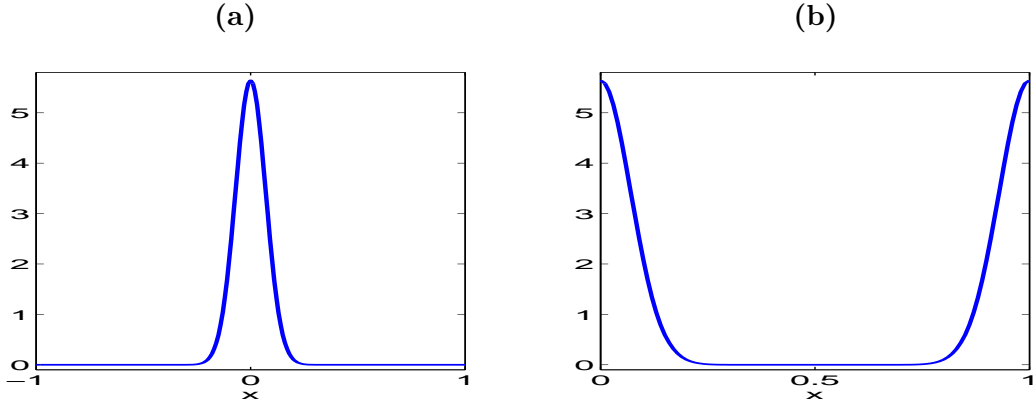}
\end{center}
\caption{The kernels $J^\infty(x)$ in (a) and $J(x)$ in (b) for the case
  of equation (\ref{eq:expk}).}
\label{fig:kernels}
\end{figure}

If we choose spatially  one-periodic
initial data $u(x,0)$, then from (\ref{eq:ut}) it is clear that
for all $x\in \mathbb{R}$ and $t\in\mathbb{R_+}$
\[
u(x,t)=u(x+1,t).
\]
Then from (\ref{eq:ut}) we have
\begin{eqnarray}\label{eq:periodic01}
u_{t}&=& \varepsilon \int_{\mathbb{R}}J^{\infty} (x-y) \left(
u(y,t) - u(x,t)\right)dy+f(u)\nonumber\\
&=&\varepsilon \sum_{r=-\infty}^{\infty}
\int_{r}^{r+1}J^{\infty}(x-y) \left(
u(y,t) - u(x,t)\right)dy+f(u)\nonumber\\
&=&\varepsilon \sum_{r=-\infty}^{\infty}
\int_{0}^{1}J^{\infty} (x-z-r) \left(
u(z+r,t) - u(x,t)\right)dz+f(u)\nonumber\\
&=&\varepsilon\int_{0}^{1}J(x-z)\left(u(z,t)-u(x,t)\right)dz+f(u),
\end{eqnarray}
where
\begin{equation}\label{eq:J}
J(x)= \sum_{r=-\infty}^{\infty} J^{\infty} (x- r)
\end{equation}
and $x\in [0, 1]$. Thus, for $1$-periodic initial data we only need
to solve the problem (\ref{eq:ut}) on the interval $\Omega=[0, 1]$
with the kernel $J(x)$. For the kernel given by (\ref{eq:expk}),
$J^\infty(x)$ and 
$J(x)$ are plotted in Figure~\ref{fig:kernels}.
%
%

\begin{lemma}\label{sym}
For $J$ defined by \eqref{eq:J} the following two properties hold.

{\em 1}. If $J^\infty(x)=J^\infty(-x)$ we have that
\[
J(x)=J(1-x).
\]
{\em 2}. $\int_0^1 J(x)\, dx =\int_{-\infty}^{\infty}
J^\infty(x)\, dx$.
\end{lemma}

Property 1. above has an important influence on the spectrum of
the matrix governing the semi-discretised version of (\ref{eq:ut}) as we
explain in the next section.
From now on we work on $[0, 1]$ and use the kernel $J$ given in (\ref{eq:J})

\section{The semi-discretised system}

We discretise in space using piecewise-constant functions~\cite{Dug}
and collocating at the uniformly spaced element mid-points, $x=
x_{j+\frac{1}{2}}$, $j=0, 1, 2, \cdots, N-1$. 
Setting $u_j= u(x_{j+\frac{1}{2}},t)$,
we have the semi-discrete approximation of (\ref{eq:ut}) given by
\begin{equation}\label{eq:utsd}
   u_t=\varepsilon  A_N u+F(u),
\end{equation}
where now $u(t) \in \R^N$, supplemented
with some initial condition $u(0)=u_0 \in \R^N$. The nonlinearity
$F: \R^N\mapsto \R^N$ is given by   $F_j(u)= f(u_j)$. 
It remains to specify the
$N \times N$ matrix $A_N$.  If we put $h = \frac{1}{N}$, its elements
are given by
\begin{equation}\label{eq:aij}
a_{j,i}=\left\{\begin{array}{ll}
        h J(|x_{j-i}|) & j\ne i \\
	h \left[J(0)-\sum_{r=1}^N \ J(|x_{j-r}|)\right ] & \ j=i.
				\end{array}
	 \right.
\end{equation}

From \lemref{sym} it follows that $A_N$ is a symmetric circulant
matrix generated by the elements $a_{1,1}$, \ldots,  $a_{1,N}$.
Hence the theory of circulant matrices can be used to characterise its
spectrum precisely. 
Let $W_k$ be the  $N$ distinct roots of $z^N-1=0$, so
$W_k=\exp{ \left(\frac{i 2 \pi k} {N}\right)}$,
for $k=0, 1, 2, ..., N-1$. Then the following theorem holds:
\begin{theorem}[\cite{K.E.Morrison}]\label{thm:eigvalues}
Let $A_N$ be the circulant matrix defined by $a_{1,1}, a_{1,2},..., a_{1,N}$.
Then $-A_N$ is diagonalisable with eigenvalues
\begin{equation}\label{eq:eigen}
\lambda_k =-[ a_{1,1}+a_{1,2}W_k+a_{1,3}W_k^2+...+a_{1,N}W_k^{N-1}],
\end{equation}
with corresponding eigenvectors $v_k =
\left(1,W_k, W_{k}^2, ..., W_{k}^{N-1}\right)^T.$
\end{theorem}

Let us see what this implies in our case for the spectrum of the
discretisation. 
\begin{lemma}\label{spec}

The following three properties hold for the spectrum of $A_N$

\textrm{1}. $\lambda_0=0$;

\textit{2}. $\lambda_{k}=\lambda_{N-k}$;

\textit{3}. Let $I_N$ be the convex hull of the set of non-zero eigenvalues
of $A_N$. As $N \rightarrow \infty$, $I_N$ converges in the Hausdorff metric
to the set
\[
I_\infty= \left[ \int_0^1 J(x)dx-\int_0^1 J(x)\exp(2\pi i x)\,
dx, \int_0^1 J(x) dx. \right].
\]
\end{lemma}
Before we prove this lemma, let us explain what it means. First of
all, we must have a zero eigenvalue with a constant
eigenvector, because, like in the case of the
Neumann Laplacian, the equation
\[
u_t= \int_{-\infty}^\infty J^\infty (x-y)(u(y,t)-u(x,t))\ dx,
\]
conserves mass.

Secondly, the pairing of the eigenvalues is simply the consequence of
the symmetry $J(x)=J(1-x)$ inherited from the evenness of the kernel
$J^\infty$. Finally, the third part of the lemma implies that as
$N\rightarrow \infty$, the spectrum accumulates at the point $\int_0^1
J(x)\, dx$. Note that in the case of $J^\infty(x)~=~\sqrt{100/\pi}~\exp(-100 x^2)$, we explicitly have
\[
I_\infty = [1-\exp(-\pi^2/100), 1] = [0.094,1].
\]

\begin{proof}
1. From (\ref{eq:aij}), putting $W_0=1$, we immediately obtain from
(\ref{eq:eigen}) that $\lambda_0=0$.

2. From part 1. of \lemref{sym} it follows that
for all $j =2,\,\ldots, \, N$,
\[
a_{1,j}=a_{1,N+2-j},
\]
so that the matrix $A_N$ is symmetric. Hence its eigenvalues $\lambda_k$ are
real. But then taking complex conjugates of $A_N v_k = \lambda_k v_k$, we
get that $A_N \overline{v_k}= \lambda_k \overline{v_k}$, or in other words
$A_n v_{N-k}= \lambda_k v_{N-k}$ and hence $\lambda_k=\lambda_{N-k}$.

3. Finally, by taking the limit as $N \rightarrow \infty$ in
(\ref{eq:eigen}) we immediately obtain that
\[
\lambda_k \rightarrow \int_0^1 J(x)\, dx - \int_0^1 J(x) \exp(2\pi i
kx) \, dx,
\]
$k=1,2,\ldots$.
\end{proof}
Our aim is to examine bifurcations in this system and, below, we
perform a numerical path--following of solution branches. Some of
these will, by symmetry, arise in pitchfork bifurcations from the
trivial solution $u=0$. Here we examine analytically the values of
$\varepsilon$ where such bifurcations may occur in the semi-discrete system 
and later we can compare to 
the numerically found values. Linearising around 
the zero solution, we have the eigenvalue problem
\begin{equation}\label{eq:eigp}
\varepsilon A_N v + \hbox{grad}\, F(0) v = \mu v,
\end{equation}
and hence bifurcations from the zero solution will only occur if
$\mu=0$, or in other words, if
\[
-A_N v = \frac{f'(0)}{\varepsilon} v.
\]
Thus, for the semi-discrete system (\ref{eq:utsd}) we can fully 
characterize the values of $\varepsilon$ where bifurcations of the 
zero solution occur, namely 
\begin{equation}\label{eq:eps}
\varepsilon_k := \frac{f'(0)}{\lambda_k}, \;
k=0,\, \ldots,\, N-1.
\end{equation}

For example, for  $N=32$,
$J^{\infty}(x)=\sqrt{\frac{100}{\pi}}\exp(-100 x^2)$ and 
$f(u)=u(1-u^2)$, we have using (\ref{eq:eigen}),
the results of \lemref{spec} and the formula (\ref{eq:eps})
that bifurcations from the zero solution are expected at the
values of $\varepsilon$ as in Table \ref{bif_egval01a:f}.
Note that for this case of $N=32$, the value of $\varepsilon_1$
agrees to 12 decimal points with the limiting
value of $\varepsilon_1$, $1/(1-\exp(-\pi^2/100))$ as $N\rightarrow \infty$,
(see part 3 of \lemref{spec}).

\begin{table}[ht]
\centering
\begin{tabular}{|c|r|c|c|r|}
\hline
$\varepsilon_1 = \varepsilon_{31}$ & $10.6403416996149 $ & &
$\varepsilon_9 = \varepsilon_{23}$ & $1.00033746722800 $ \\
\hline
$\varepsilon_2 = \varepsilon_{30}$ & $3.06584313146254 $ & &
$\varepsilon_{10} = \varepsilon_{22}$ & $1.00005172586163 $ \\
\hline
$\varepsilon_3 = \varepsilon_{29}$ & $1.69885748860222 $ & &
$\varepsilon_{11} = \varepsilon_{21}$ & $1.00000650971543$ \\
\hline
$\varepsilon_4 = \varepsilon_{28}$ & $1.25968856776757 $ & &
$\varepsilon_{12} = \varepsilon_{20}$ & $1.00000067252291$ \\
\hline
$\varepsilon_5 = \varepsilon_{27}$ & $1.09266327932320 $ & &
$\varepsilon_{13} = \varepsilon_{19}$ & $1.00000005703325$ \\
\hline
$\varepsilon_6 = \varepsilon_{26}$ & $1.02948119722480 $ & &
$\varepsilon_{14} = \varepsilon_{18}$ & $1.00000000397031$ \\
\hline
$\varepsilon_7 = \varepsilon_{25}$ & $1.00800141737791 $ & &
$\varepsilon_{15} = \varepsilon_{17}$ & $1.00000000022729$ \\
\hline
$\varepsilon_8 = \varepsilon_{24}$ & $1.00180943793728 $ & &
$\varepsilon_{16}$ & $1.00000000002128$ \\
\hline
\end{tabular}
\label{bif_egval01a:f}
\caption{For $N=32$ Bifurcation values in terms of $\varepsilon$ of
  the zero solution.}
\vspace*{0.2in}
\end{table}

Let us examine the eigenvectors of $-A_N$ in some more detail.
Since both $v_k$ and $v_{N-k}$ are eigenvectors, we
immediately have that $\hbox{Re}\,(v_k)$ and $\hbox{Im}\,(v_k)$ are
eigenvectors. Define the cyclic shift $\sigma$ on $u = (u_1,\,\ldots,\,u_N) \in
\R^N$ by
\[
\sigma(u) = (u_N,\, u_1,\, \ldots,\, u_{N-1}),
\]
then we have
\begin{lemma} \label{eigenv}
If $v$ is a real eigenvector of $-A_N$ corresponding to a double
eigenvalue $\lambda$, then so  is $\sigma(v)$.
\end{lemma}
This follows since if $v$ is an eigenvector, then so is
$e^{2i\pi/N}v$.

{\bf Remark.} In the above argument, we can pass to the limit as $N
\rightarrow \infty$ and arrive at the somewhat startling conclusion that
$\cos(2\pi kx)$ and all their translates are eigenfunctions of
$-A=-(\int_0^1 \, J(x-y) (u(y)-u(x))\, dy)$ {\em no matter what the
kernel $J(x)$ is\/} as long as it has the right symmetry property. Of course,
cosines are also the eigenfunctions of the Neumann Laplacian. It is very
pleasing to obtain such a result via a semi-discretisation.

Finally we note that fixed points of the semi-discrete problem satisfy
\begin{equation}
\label{eq:sd}
0=\varepsilon A_N u + F(u).
\end{equation}
Thus at $\varepsilon = 0$ stable solutions are given by
\begin{equation}\label{f:definidata01}
u = \left\{
     \begin{array}{rr}
  1 & \mbox{$x \in X$},\\
-1 &\mbox{$x \in Y$}
\end{array}
\right.
\end{equation}
where $X \cup Y = [0, 1]$.
Unstable solutions at $\varepsilon = 0$ are given by
\begin{equation}\label{f:definidata01a}
u = \left\{
     \begin{array}{rr}
  1 & \mbox{$x \in X$},\\
-1&\mbox{$x \in Y$},\\
0 &\mbox{$x \in Z$}
\end{array}
\right.
\end{equation}
where $X \cup Y\cup Z = [0, 1]$ with some nonempty $Z$.

We use this to define solutions with different numbers of interfaces.
When $X=[0, \alpha)$ and $Y=[\alpha, 1]$, $0 < \alpha < 1$, we call 
$u$ a one-interface solution of (\ref{eq:sd}) if for $\varepsilon=0$ 
for some $n \in
[0,N-1]$, $\sigma^n(u)=a_1$ on $X=[0,
\alpha)$ and $\sigma^n(u)=a_2$ on $Y=[\alpha, 1]$, $a_1\, a_2 \in \{-1, 0,1\}$,
$a_1 \neq a_2$. That is loosely speaking we have at $\varepsilon=0$
one jump in the solution upto cyclic shift.
Two-interface, three-interface solutions, etc., are defined
similarly. 
Thus, for example, the branch of solutions corresponding to orbit $A$ in
Table \ref{tab:2} are of one-interface and those corresponding to $E$
are of three-interface.


\section{Results}\label{results_bif:f}
%
We take for our computations the kernel function  $$J^{\infty}
(x)=\sqrt{\frac{100}{\pi}}e^{-100x^2},$$
with $f(u) = u(1-u^2)$ and vary the parameter $\varepsilon$.
For small values of $N$ it is possible to enumerate all possible
solutions of the semi-discrete system (\ref{eq:sd}) with $\varepsilon=0$
and to analyse their continuation to $\varepsilon>0$ using the theory 
of bifurcation with symmetry.
This we do below for $N=4$ and these analytic results were used 
to check the validity of our numerics.

We implemented in \mat a standard pseudo arc--length continuation algorithm
with step size control as described in~\cite{Govaerts, Z.mei, W.Gov00}
for the discrete problem (\ref{eq:sd}).
Since $A_N$ is a circulant matrix, we take advantage of reducing storage 
costs as the full information of $A_N$ can be obtained storing one row or 
column only, see \cite{SKB} and references therein. Furthermore the use 
of the FFT for each matrix vector multiplication reduces the
computational cost. 
%
We detect bifurcation points by observing where eigenvalues of
the Jacobian $\mathcal{J} = D_u F$ of the nonlinear system
$F(u,\beta)=0$ cross the imaginary axis and perform branch
switching at those points by perturbing in the direction of the
associated eigenvector.

The arc-length $\ell$ of $u(x)\in C^1(\mathbb{R})$ is defined in the
standard way
\[
\ell=\int_{\Omega} \sqrt{1+\left (\frac{du}{dx}\right )^2} dx,
\]
and we approximate the arc--length of $u(x)$ with the mid-point rule and
using the standard forward difference approximation for the derivative.
With a uniform discretization we get
\begin{equation}\label{f:arc_lintthh}
\ell \approx \ell_h = \sum_{j=0}^{N-1} \sqrt{h^2 + \left (u_{j+1}-u_{j}\right )^2}
\end{equation}
where $h = x_{j+1}-x_j$ and $u_j\approx(x_j)$.
Note that although $\ell$ only makes sense for $u\in C^1$ however we can
evaluate $\ell_h$ even when $u$ is discontinuous at grid points.

Then, for $N=32$ we compute the bifurcation diagram numerically and gain
insight into the structure of the bifurcation diagram of the original
continuous problem.

Finally, we examine the large $N$ limit and formulate the results of the
numerics as two conjectures concerning the interplay of continuity and
stability and the behaviour of saddle-node bifurcations as $\alpha
\rightarrow 1/2$.


For the continuous system, the symmetry group is $O(2) \times \Z_2$,
and so for a finite number of nodes $N$, we use $\Gamma_N = D_N
\times \Z_2$ equivariance structure \cite{GSS,Hoyle}.

\subsection{The $N=4$ case}

If $N=4$, there are a total of $81$ possible steady states at
$\varepsilon=0$, $16$ of them stable. The group $\Gamma_4$ is
generated by the shift 
$p$, the flip $f$ and the reversal $m$. In other words, if $v=(v_1,\, v_2,\,
v_3,\; v_4)$, we have  that
\[
\begin{eqalign}
&p(v)=(v_2,\, v_3,\, v_4,\, v_1);\\
&f(v)=(v_4,\, v_3,\, v_2,\, v_1);\\
&m(v)=(-v_1,\, -v_2,\, -v_3,\, -v_4);
\end{eqalign}
\]
Inverses of the nonzero eigenvalues of the $4 \times 4$
matrix $A_4$ are $\{91.82, 183.63, 183.63 \}$, so we expect primary branches
to bifurcate from the zero solution at those values of $\varepsilon$. Note
that all primary branches have zero mean, but the converse
is not true.

Since here we know all the solutions at $\varepsilon=0$ and their stability,
and since symmetry properties are conserved on primary branches, we can cut
down the work considerably by looking only at orbits of solutions under
$\Gamma_4$. In the table \ref{tab:2}, we collect all the orbits, their lengths and
the corresponding isotropy subgroups $\Sigma_x$. There, $\ip{f}$ stands for
the group generated by $f\in D_4 \times \Z_2$. 

\begin{table}
\begin{tabular}{|c|r|r|c|}
\hline
Name & Orbit        & length & $\Sigma_x$\\ \hline
     & (0,0,0,0)    &   1    & $D_4 \times \Z_2$ \\
     & (1,1,1,1)    &   2    &  $D_4$\\ \hline
$A$    & (-1,1,1,-1)  &   4    &   $\ip{f,p^2 m}$\\
$B$    & (1,0,-1,0)   &   4    & $\ip{mp^2,fp}$\\
$C$    & (0,0,1,1)    &    8   & $\ip{f}$\\
$D$    & (0,-1,1,1)   &   16   & $\ip{I}$ \\ \hline
\end{tabular}
\begin{tabular}{|c|r|r|c|}
\hline
Name & Orbit        & length & $\Sigma_x$\\ \hline
$E$    & (-1,1,-1,1)  &   2    & $\ip{pm,fp}$\\
$F$    & (0,1,0,0)    &    8   & $\ip{fp}$\\
$G$    & (1,0,0,-1)   &   8    & $\ip{mf}$\\ 
$H$    & (0,1,0,1)    &    4   & $\ip{p^2,fp}$\\
$I$    & (0,1,-1,1)   &    8   & $\ip{fp}$\\\hline
$J$    & (0,1,1,1)    &   8    & $\ip{fp}$ \\
$K$    & (-1,1,1,1)   &   8    & $\ip{fp}$ \\ \hline
\end{tabular}
\caption{Steady states for $N=4$, the length of the orbits and
  isotropy subgroups $\Sigma_x$. We have separated the solutions 
  into the homogeneous
  states, those connected with the first and second bifurcation of
  $u=0$ and the two solutions connected by a saddle-node $J$ and $K$.
  See also Figure \ref{fig:91and183}.}
\label{tab:2}
\end{table}

Now we can immediately draw the bifurcation diagram using the following
three rules \cite{GSS,Hoyle}. First a bifurcating branch must have the
isotropy subgroup which is a subgroup of the isotropy subgroup of the
primary branch; secondly dimensions of unstable manifolds have to
match at a bifurcation point to 
satisfy the principle of exchange of stability, and thirdly at
$\varepsilon =0$; the number of nodal domains must increase from one
bifurcation point to the next.

\begin{figure}
\begin{center}
{\bf (a) \hspace{0.48\textwidth} (b)}
\includegraphics[width=0.48\textwidth,height=6.5cm]{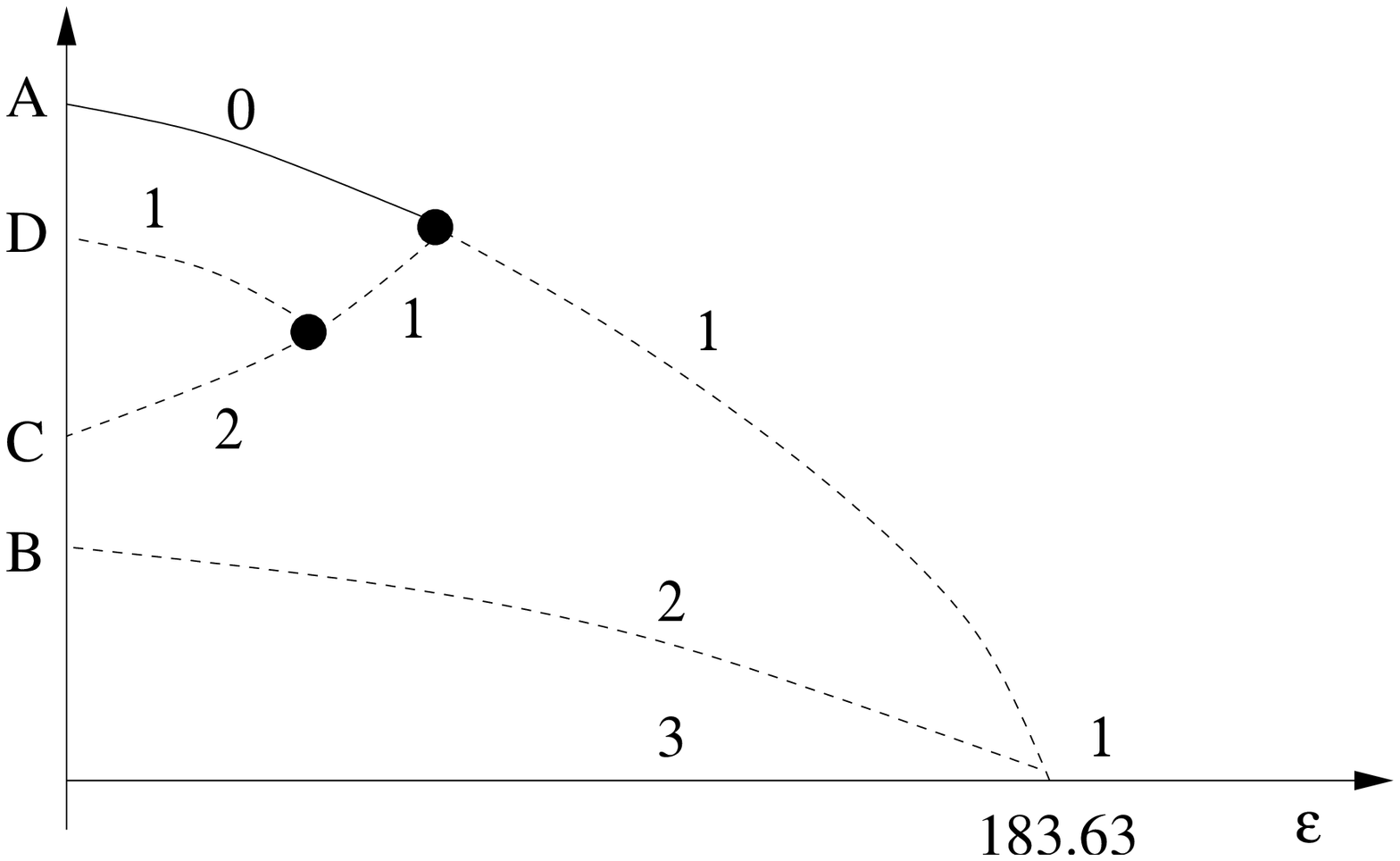}
\includegraphics[width=0.48\textwidth,height=6.5cm]{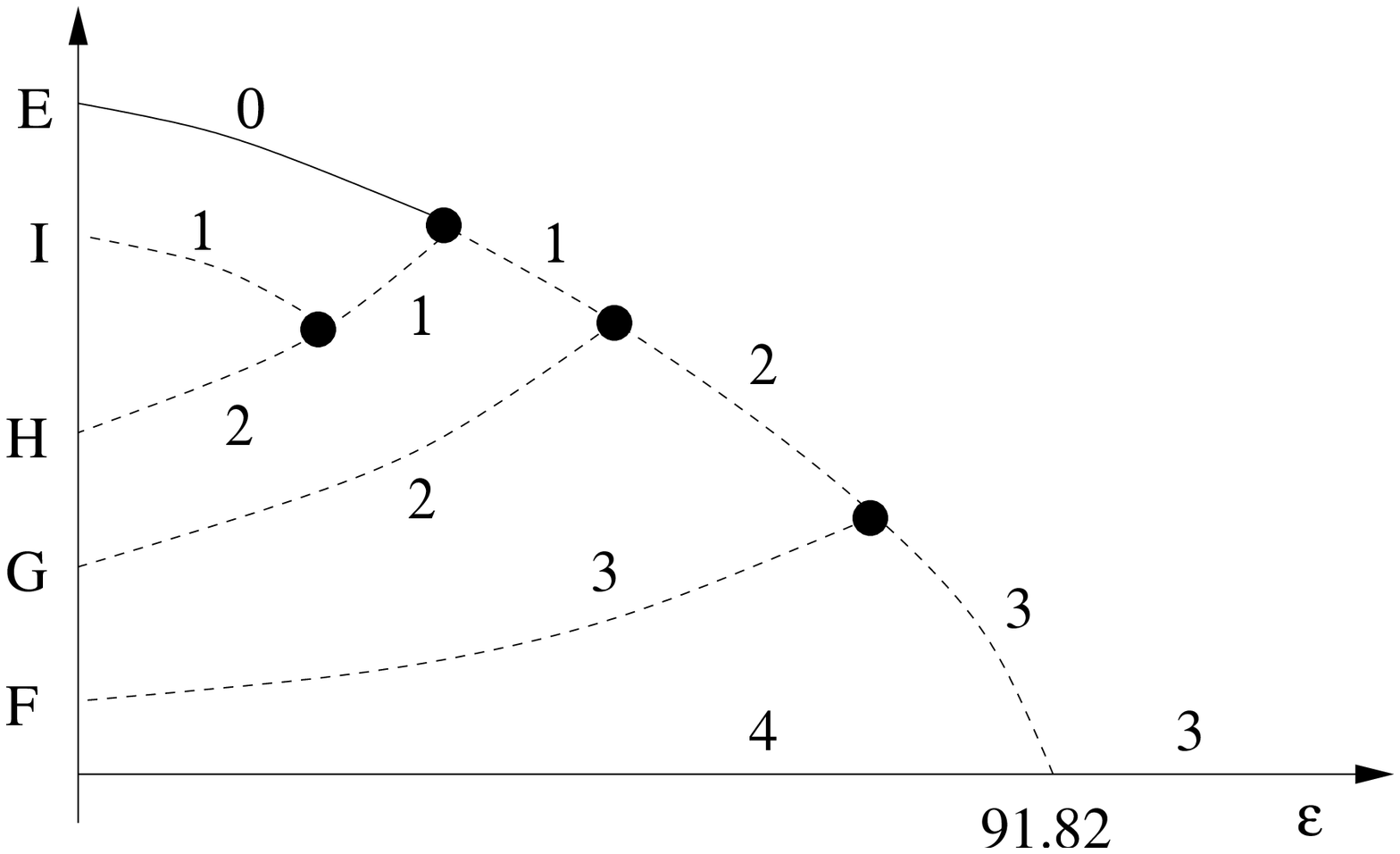}
\end{center}
\caption{Bifurcations from the zero solution from the first (a) and
  second (b) bifurcation points. Solid lines represent stable and
  broken lines unstable unstable solutions respectively. We indicate
  above each branch the dimension of the unstable manifold. We do not
  show the homogeneous solutions or the orbits $J$ and $K$ (which are connected
  by a saddle-node bifurcation).}
\label{fig:91and183}
\end{figure}

With these rules there is only one way to construct the bifurcation
diagram; see Figure \ref{fig:91and183} (a) and (b), 
where the
$y$-axis is not to any scale, and is only intended to make clear the end-points
of various branches at $\varepsilon=0$. These figures show the bifurcation
structure arising from bifurcations of the zero solution. 

\begin{figure}
  \psfrag{A}{A}
  \psfrag{B}{B}
  \psfrag{C}{C}
  \psfrag{D}{D}
  \psfrag{E}{E}
  \psfrag{F}{F}
  \psfrag{G}{G}
  \psfrag{H}{H}
  \psfrag{I}{I}
  \psfrag{J}{J}
  \psfrag{K}{K}
  \psfrag{epsilon}{$\epsilon$}
\begin{center}
{\bf (a) \hspace{0.49\textwidth} (b)}
\includegraphics[width=0.49\textwidth,height=6.5cm]{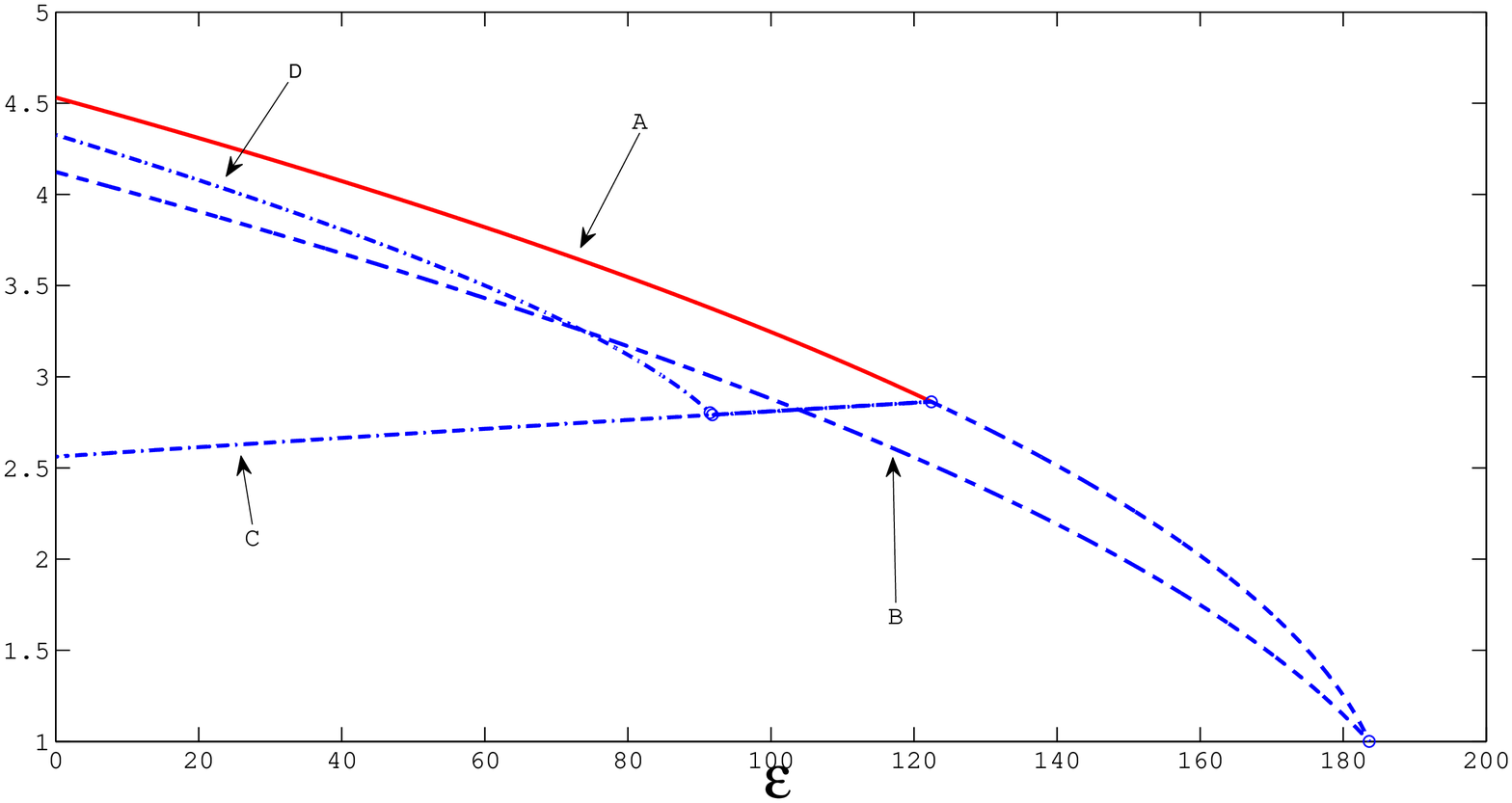}
\includegraphics[width=0.49\textwidth,height=6.5cm]{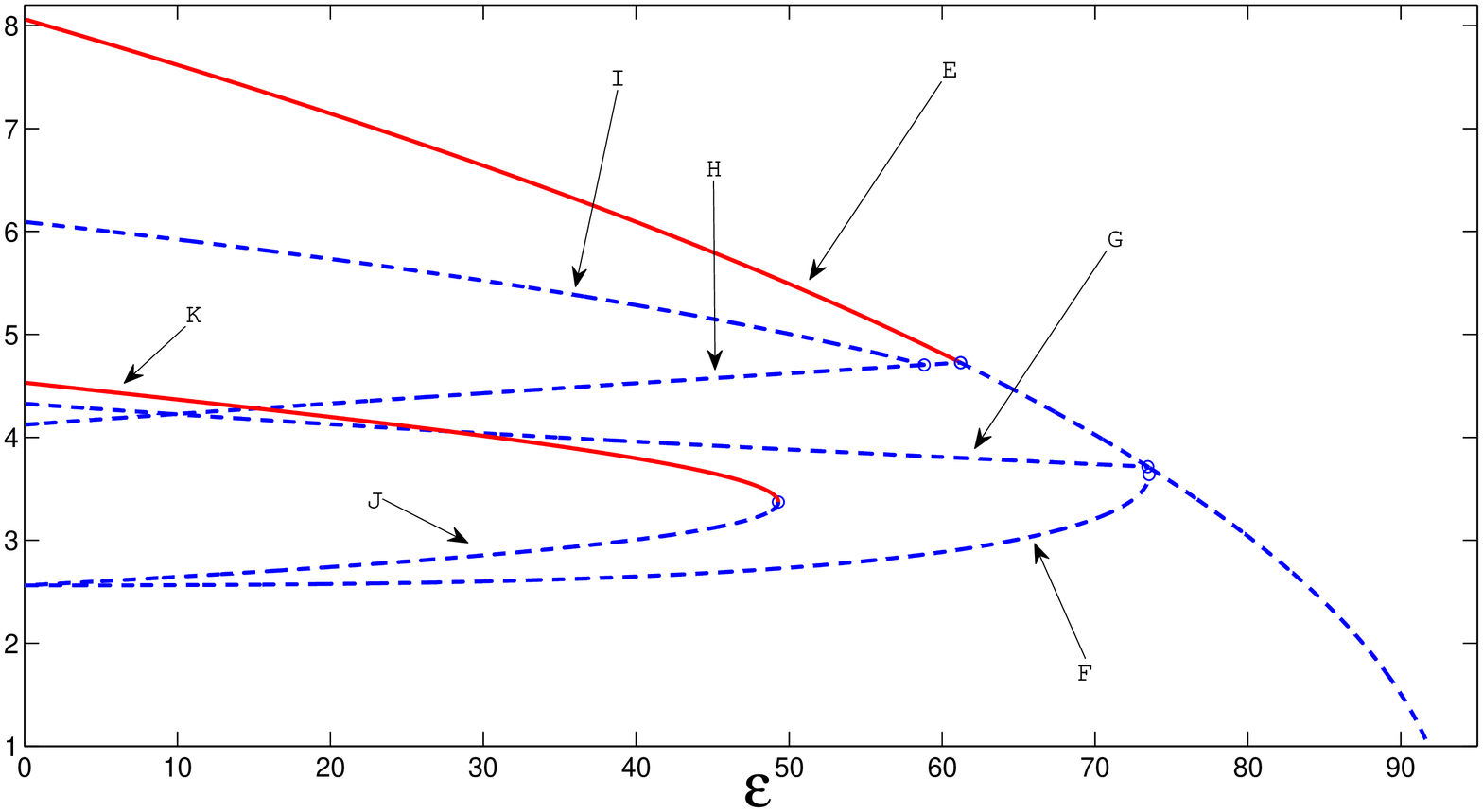}
\end{center}
\caption{Numerically computed bifurcation diagram from the zero solution from the first (a) and
  second (b) bifurcation points. Solid lines represent stable and
  broken lines unstable unstable solutions respectively. We indicate
  above each branch the dimension of the unstable manifold. 
  We show the homogeneous solutions or the orbits $J$ and $K$ connected
  by a saddle-node bifurcation on (b). Compare to the theoretical
  prediction in Figure \ref{fig:91and183}.}
\label{fig:Num91and183}
\end{figure}

We would like to make the following observations. The stable non-zero-mean
branches corresponding to the orbit $K$ have to arise through a saddle-node
bifurcation. Numerically, this happens at a value of epsilon $\approx
49.294$ that is smaller than the value $\varepsilon_4=122.432$ at which
the branches of the orbit $A$ become stable, see Figure
\ref{fig:Num91and183} which shows the numerically computed
diagram. We will  see the equivalents of these statements in higher
dimensional discretisations.  

Finally, we did not perform a Liapunov--Schmidt calculation to determine the
order of bifurcations at the double eigenvalue point $\varepsilon=183.63$,
but the opposite assignment of stabilities cannot be reconciled with the above
rules of bifurcation.

\subsection{The case of $N=32$}
\begin{figure}[here]
  \psfrag{si}{seven-interface}
  \psfrag{fi}{five-interface}
  \psfrag{ti}{\hspace{3ex}three-interface}
  \psfrag{oi}{one-interface}
\begin{center}
{\bf (a) \hspace{0.49\textwidth} (b)}\\
\includegraphics[width=0.49\textwidth]{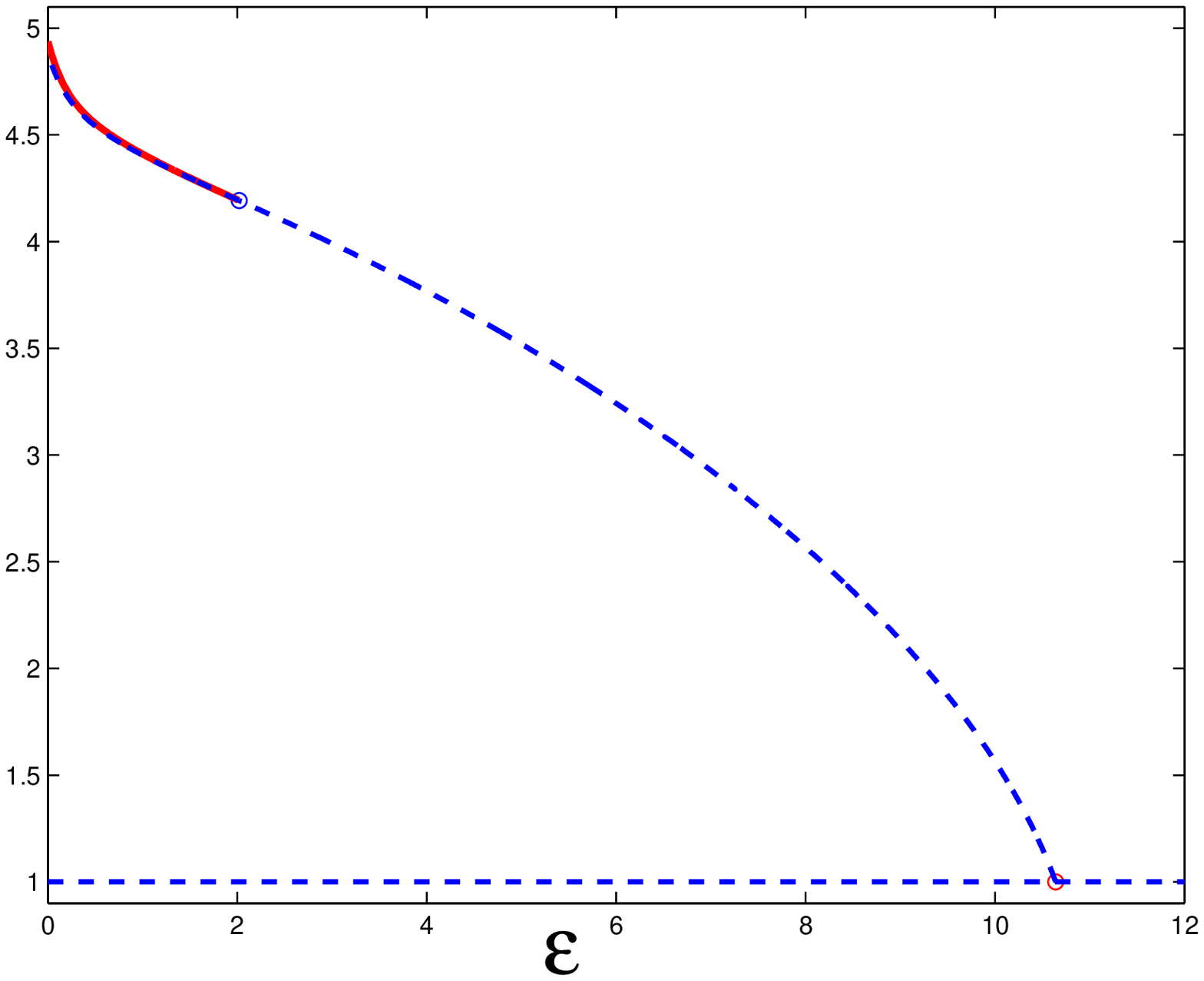}
\includegraphics[width=0.49\textwidth,height=0.27\textheight]{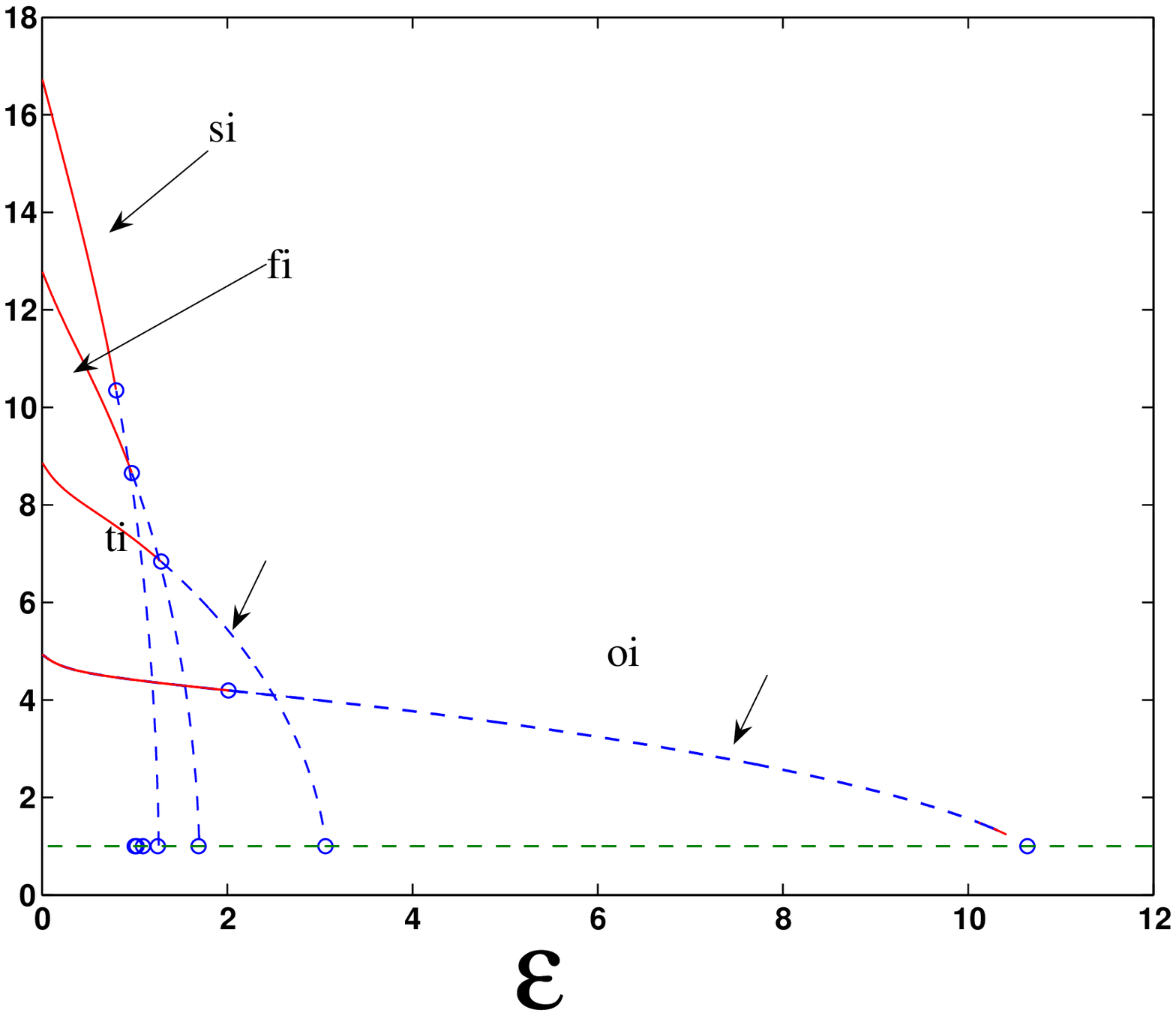}
{\bf (c) \hspace{0.49\textwidth} (d)}\\
\includegraphics[width=0.49\textwidth]{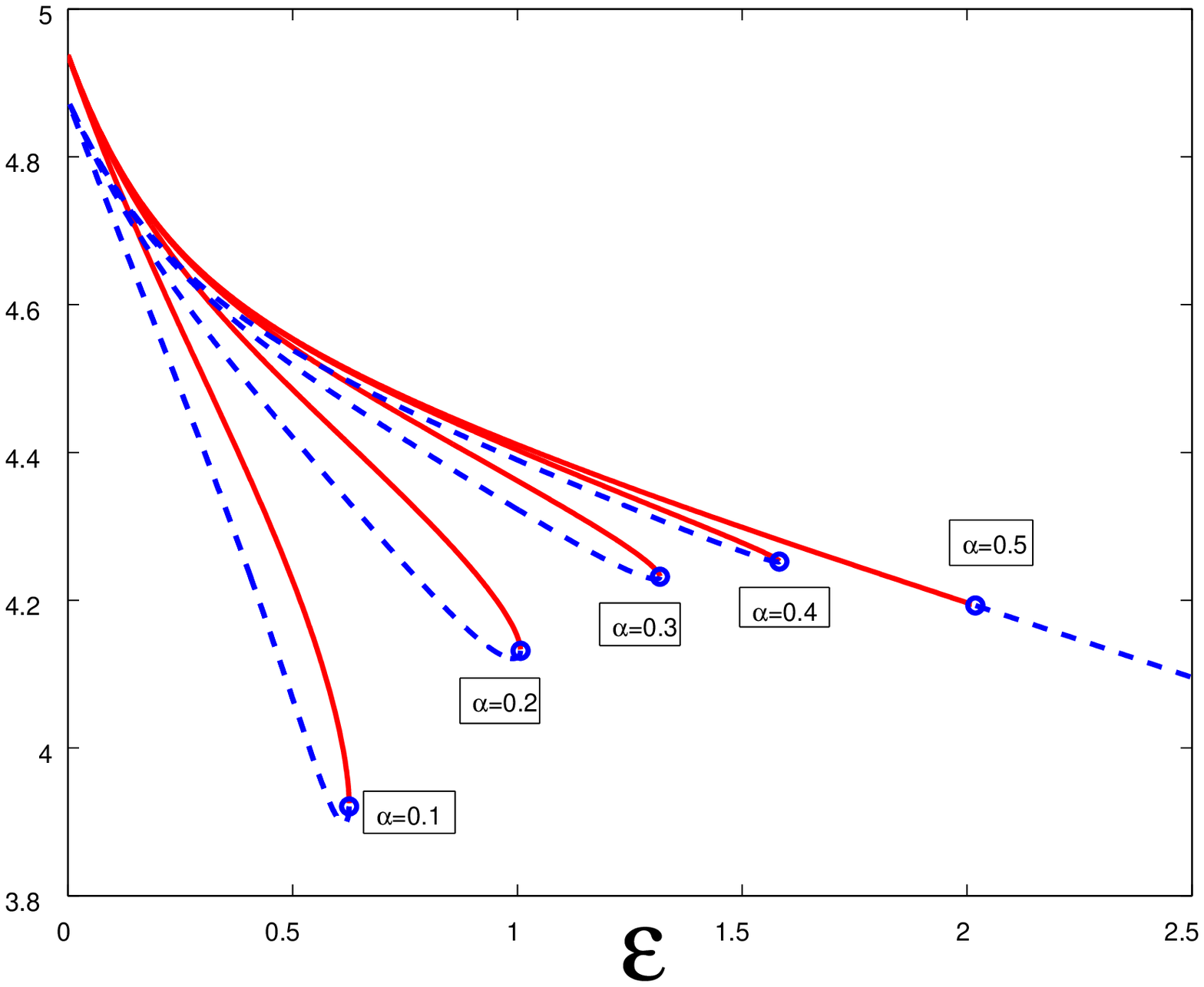}
\includegraphics[width=0.49\textwidth]{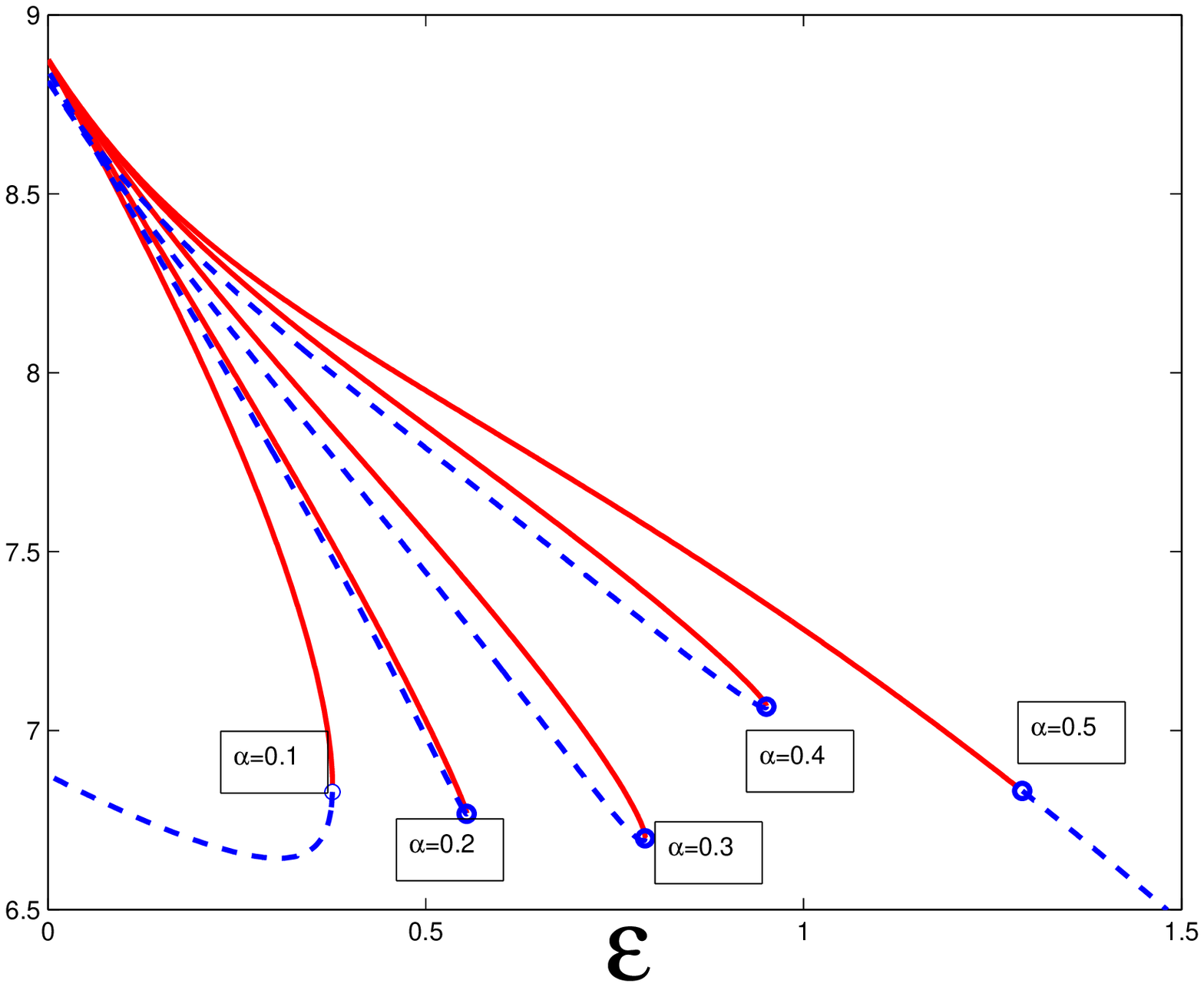}
\caption{(a) Continuation of the zero-mean one-interface ($\alpha=0.5$)
solution branch. When $\varepsilon \approx 2.012$ shown by 'o' there
is a stabilizing pitchfork bifurcation.
(b) Stabilization of one-, three-, five-, and seven-interface solutions.
(c) This is a blow up of the bifurcation diagram around the
stabilizing bifurcation at $\varepsilon \approx 2.012$ shown by 'o' for the
$\alpha=0.5$ curve. Here stable solutions are continued from
$\varepsilon=0$ with different ratios ($\alpha$) of $-1$ and $1$ values of $\alpha$.
(d) A similar structure is observed starting from $\varepsilon=0$ with
two--interface solutions with different ratios ($\alpha$) of $-1$ and $1$.
}
\label{fig:n32}
\end{center}
\end{figure}

Though an analysis similar to that in the case of $N=4$ can be attempted
here, the numbers of orbits are astronomical, and we rely on our numerical
continuation method, the results of which match exactly the predictions of
the analysis in the case $N=4$.
In Figure \ref{fig:n32} we plot in (a)--(d) sample solution branches of the
bifurcation diagram with $N=32$. If we start with a large value of
$\varepsilon$ we see in (a) and (b) the first bifurcation arises at
$\varepsilon \approx 10.64$ as predicted by the theory in Table
\ref{bif_egval01a:f}. In (a) we show the continuation of the zero mean
one-interface which undergoes a pitchfork bifurcation at $\epsilon
\approx 2.012$. In (b) we have
plotted the one-, three-, five- and seven interfaces and their
stabilization.
In (c) we show details of the bifurcation structure close to the
pitchfork at $\varepsilon \approx 2.012$ (note for clarity one branch
of the pitchfork seen in (a) is not plotted). 
As $\alpha \rightarrow 0.5$, the saddle-node bifurcation points
converge to $\varepsilon_{3,2}$. This structure is repeated
for the other $n$-interface solutions and is illustrated in (d) for
the three-interfaces solutions.
Here we see that the zero-mean one-interface solution branches stabilize at
$\varepsilon_{3,2} \approx 2.012$. Below we will formulate a conjecture
concerning the limiting value which we call $\varepsilon_1^0$ at
which the one-interface branch with zero-mean stabilizes as $N
\rightarrow \infty$.

\section{The limiting problem and conclusions}

It is not hard to prove (see for example \cite{C.Fife}) that if
$\varepsilon>1$, steady state solutions of (\ref{eq:ss}) are continuous,
since the function $-\varepsilon u +f(u)$ is monotone.  Hence it is
interesting to understand when the solutions lose continuity (certainly, for
$\varepsilon=0$ there are no non-constant continuous solutions).

The non-trivial stable one-interface zero-mean solution branches
($\alpha=0.5$) that originate at $\varepsilon=0$ were investigated in detail
as we change $N$. If we define $M$ by
$$M := \max_j \left|\frac{u_{j+1}-u_j}{h}\right|$$
then for a $C^1$ function this converges to $\max_{x\in[0,1]}
\left|u_x\right|$ and so we can identify where the solution is continuous.

\begin{figure}[here]
\begin{center}
{\bf (a) \hspace{5.5cm} (b)}\\
\includegraphics[width=0.49\textwidth,height=7.5cm]{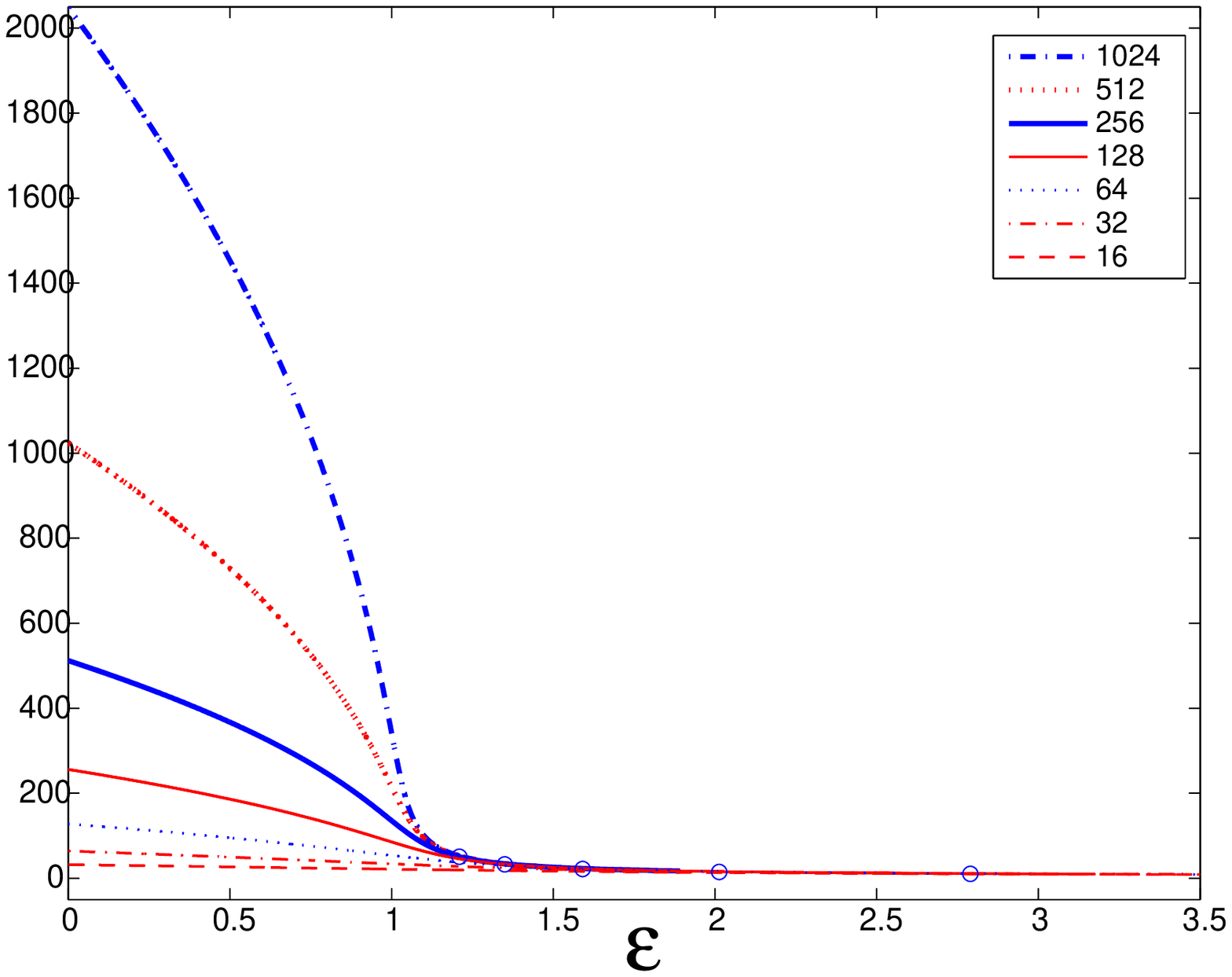}
\includegraphics[width=0.49\textwidth,height=7.5cm]{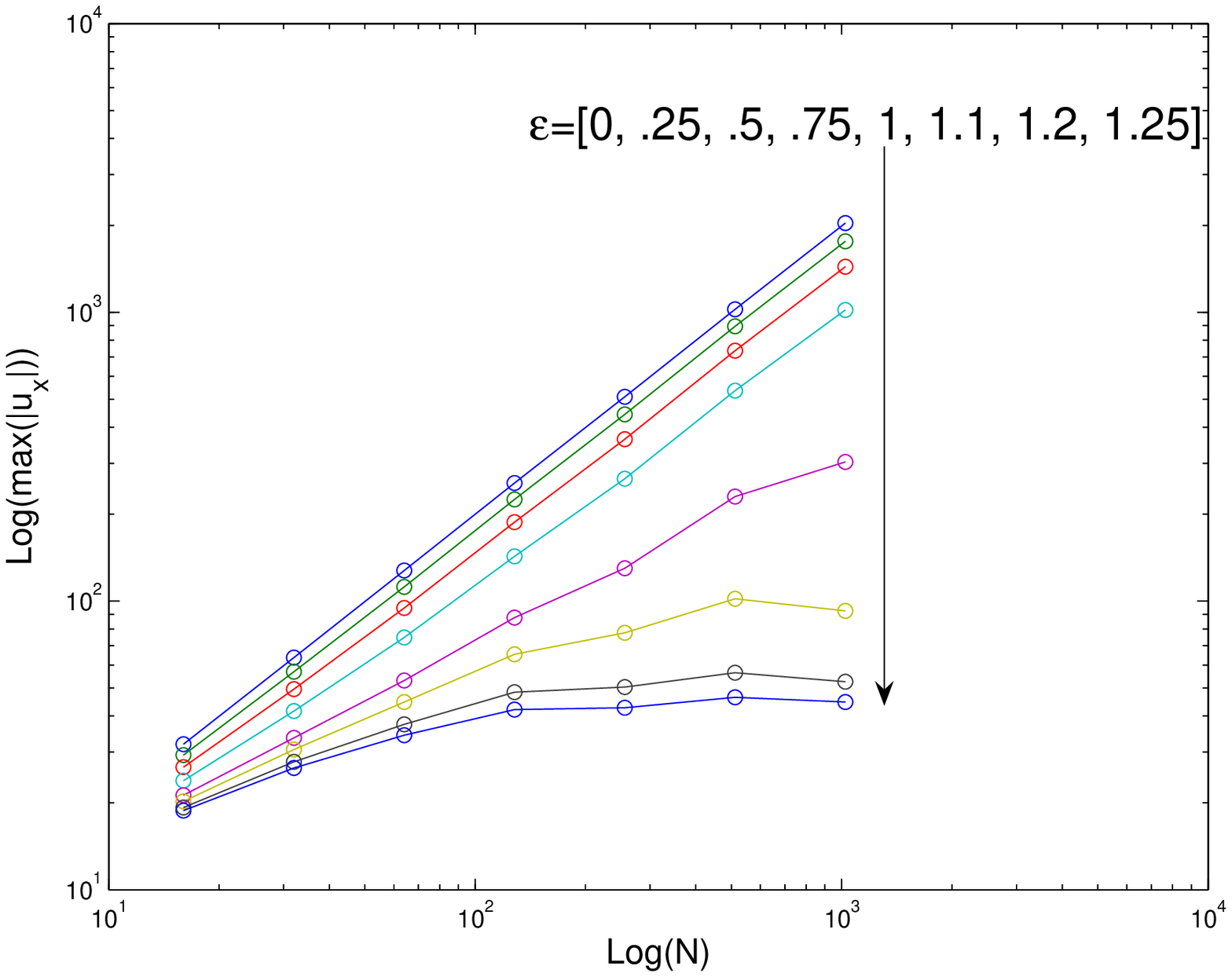}
\caption{In (a) we show the maximum of $|u_x|$ in $[0, 1]$
of the solution $u (x)$ for one-interface initial data with $\pm 1$
and $\alpha=0.5$. 
In (b) we see $\max_x(u_x)$ for different system sizes for some particular
$\varepsilon$'s with system size as the x-axis with a clear change in
behaviour at $\varepsilon=1$.}
\label{fig:ux}
\end{center}
\end{figure}

Figure \ref{fig:ux} plots in (a) $M$ against $\varepsilon$ along a branch of
one-interface zero-mean solutions for $N=2^p$,
$p=4,5,6,7,8,9,10$. 
If we let $\varepsilon_1^d$ be the value of $\varepsilon$ at
which this branch of solutions becomes discontinuous, 
then this figure suggests $\varepsilon_1^d=1$. 
This is supported in (b)
which shows for different $\varepsilon$ convergence of the derivative
$M$ with $N$ on a $\log\log$ scale.

\begin{figure}[here]
  \psfrag{bifurcation}{$\varepsilon_m^0$}
  \psfrag{bifurcation1m}{$\varepsilon_m^0-z$}
\begin{center}
{\bf (a) \hspace{5.5cm} (b)}\\
\includegraphics[width=0.49\textwidth,height=7.5cm]{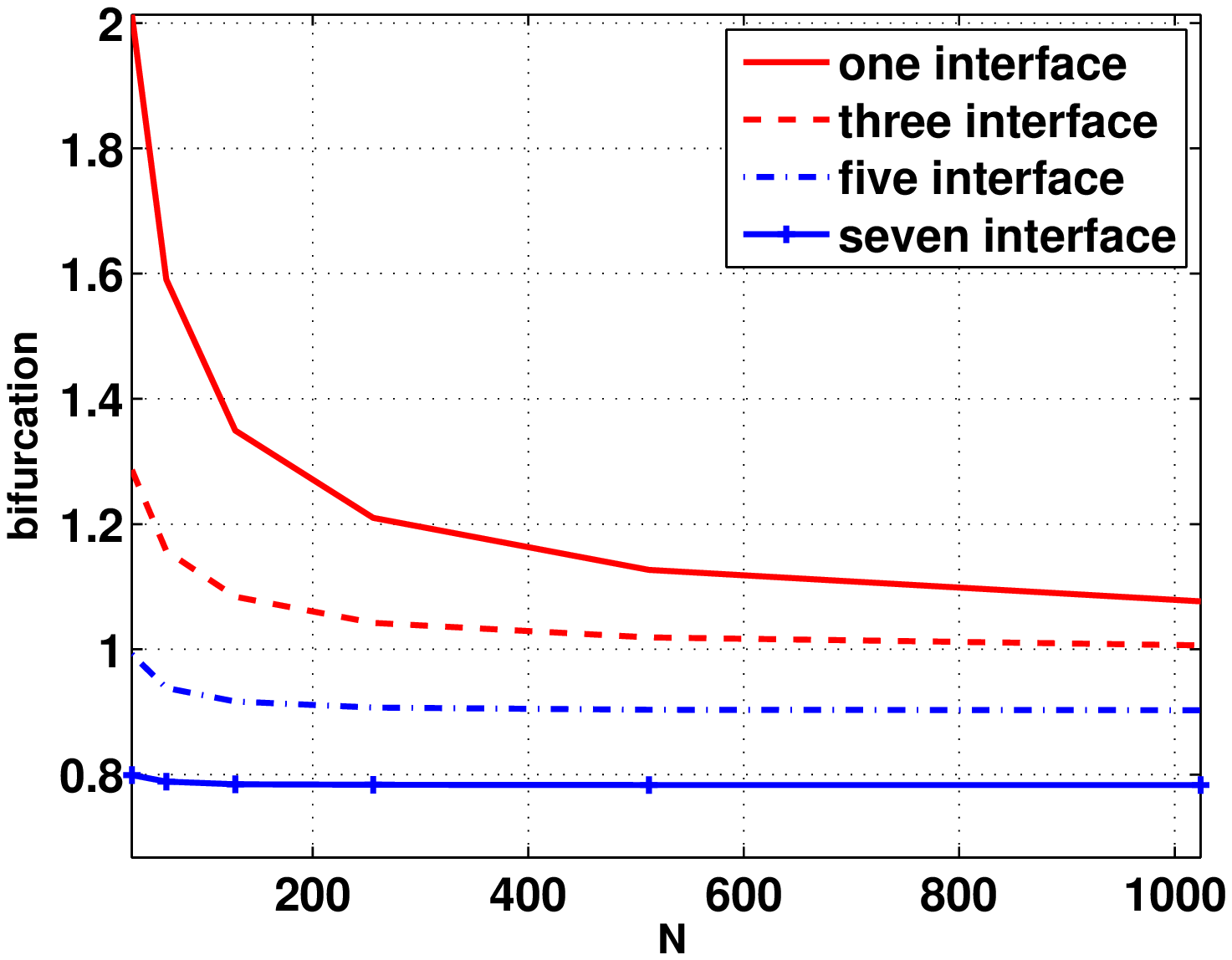}
\includegraphics[width=0.49\textwidth,height=7.5cm]{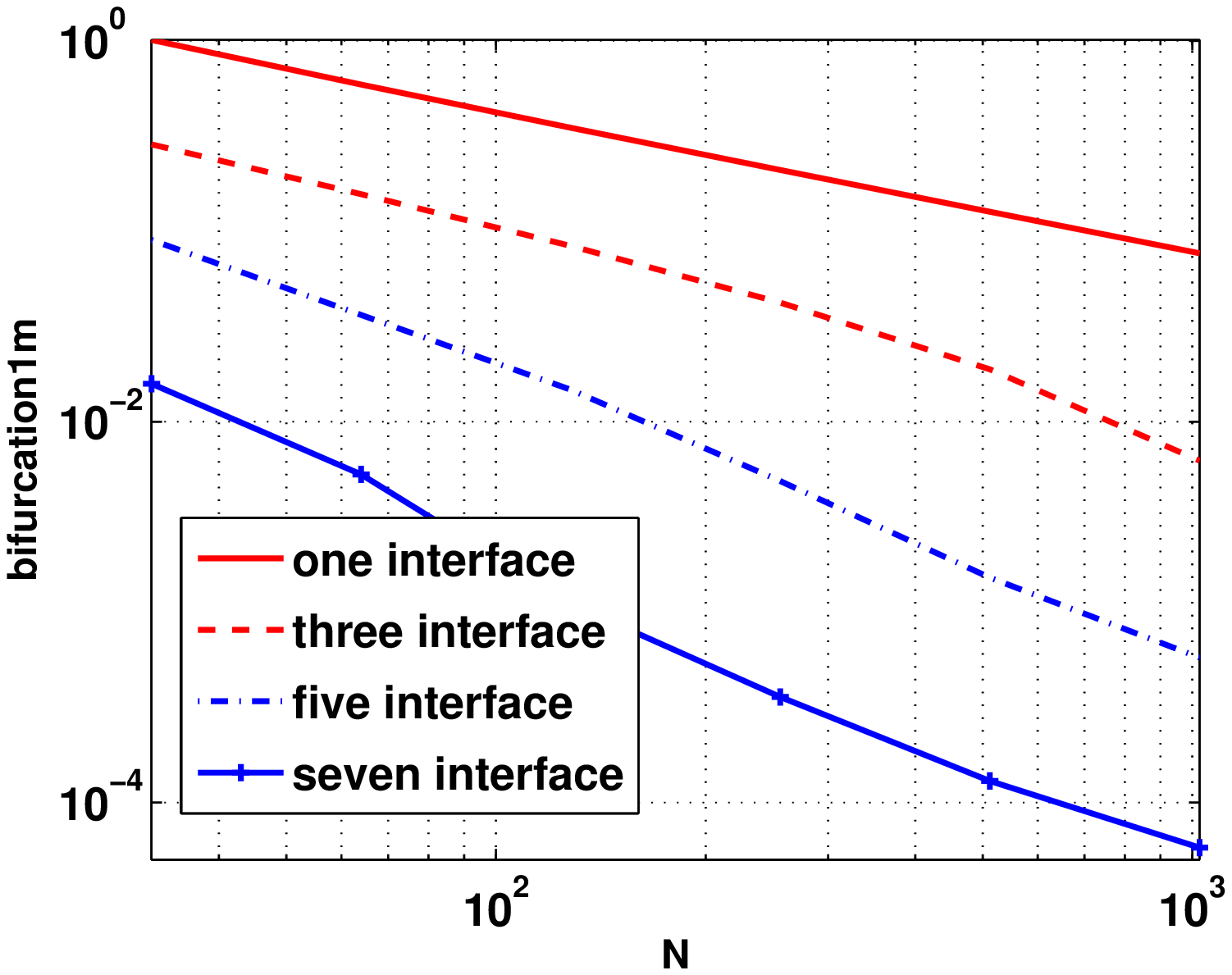}
\caption{(a) Pitchfork bifurcations of the first $(2m-1)$--interface
  solutions as $N$ increases. (b) A log-log plot of  $ \varepsilon -
  z$ versus $N$ showing convergence to $z=1, 1, 0.9, 0.7839$ for the
  one, three, five and seven--interface solutions.
}
\label{fig:bif}
\end{center}
\end{figure}

Furthermore the loss of continuity appears to coincide with the loss of
stability. In Figure \ref{fig:bif} we show numerically that the
bifurcation values converge to $\varepsilon_1^0= \varepsilon_2^0=\varepsilon_3^0=1$ 
as $N \rightarrow \infty$. 

Now we can collect our observations and form two conjectures.  
First we consider the zero-mean interface branches.
Let $\varepsilon_m^0$ be the value at which the zero-mean $2m-1$-interface
solution becomes stable. 
Let $\varepsilon_m^d$ be the value of $\varepsilon$ at
which this branch becomes discontinuous.  Then we have

{\bf Conjecture 1:} $\varepsilon_m^0= \varepsilon_m^d$.

We can prove a very weak form of this conjecture for $m=1$. From the results of
\cite{C.Fife} it follows that discontinuous stationary solutions will exist
for any $\epsilon$ such that the function 
\[
g(u) := -\varepsilon u \int_0^1 J(s) \, ds + f(u)
\]
is non-monotone. On the other hand, from Theorem 2.1 of \cite{CR02} it follows
if $g(u)$ is monotone, there are no nonconstant minimizers of the energy
functional (\ref{fe:f}). Hence we have $\varepsilon_1^s \leq  \varepsilon_1^d$.
However we do not have the inequality the other way.

We now consider the saddle-node bifurcation of the non--zero mean
interface solutions. 
Now, let $u_s$ be a branch of $2m-1$-interface stable solutions of
(\ref{eq:ss}) with mean $s$, and let $\varepsilon^{b,s}_m$ be the value of
$\varepsilon$ at which the saddle-node bifurcation giving rise to the branch
occurs. Then we have 

{\bf Conjecture 2:} 
${\displaystyle \lim_{s \rightarrow 0} \varepsilon_m^{b,s} =
\varepsilon_m^0}$. 

\begin{figure}[here]
\begin{center}
  \psfrag{LPs}{$\varepsilon_m^{b,s}$}
  \psfrag{e10e1d}{$\varepsilon_1^{0}=\varepsilon_1^{d}$}
  {\bf (a) \hspace{5.5cm} (b)}\\
  \includegraphics[width=0.49\textwidth,height=7.5cm]{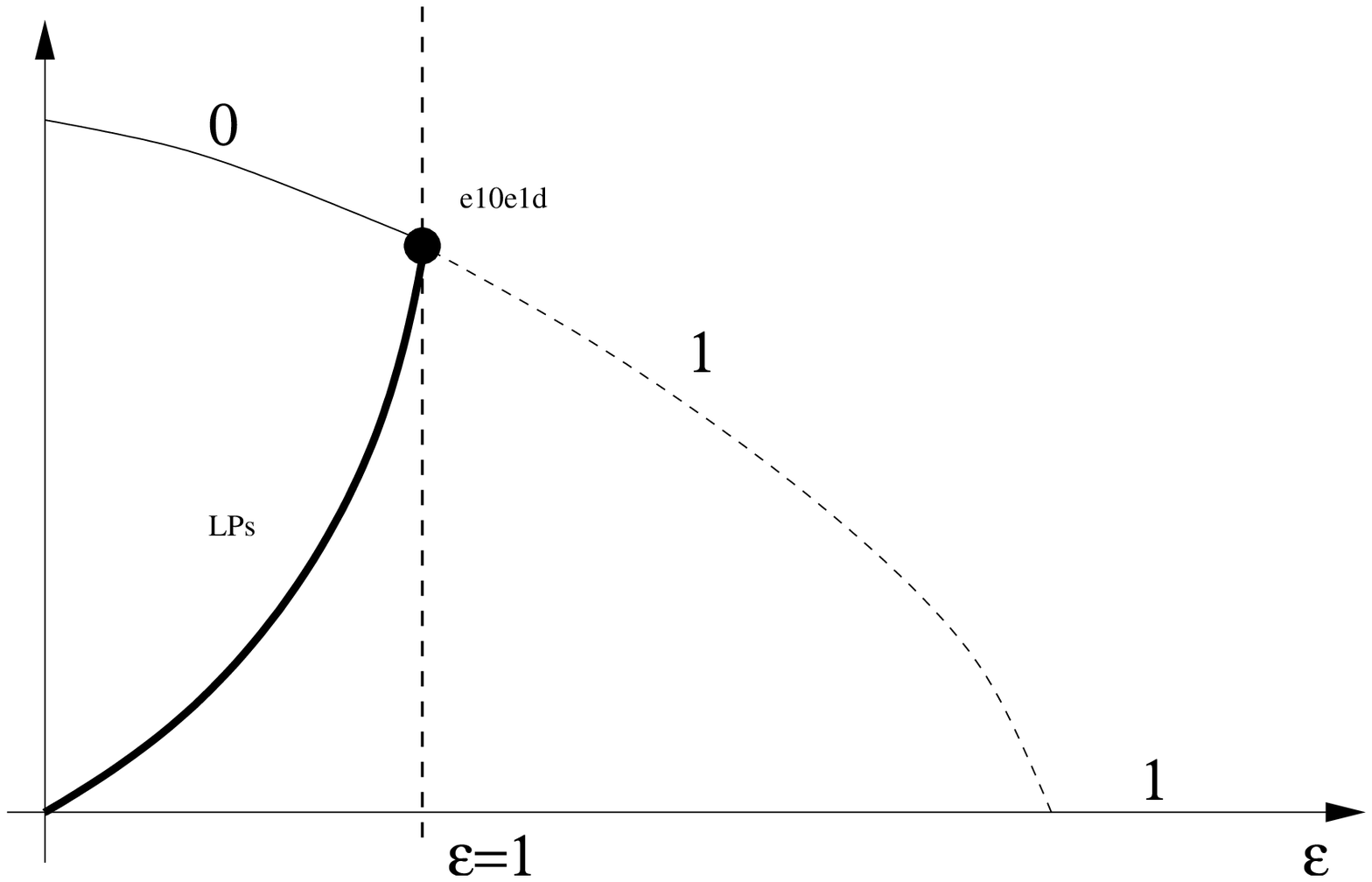}
  \includegraphics[width=0.49\textwidth,height=7.5cm]{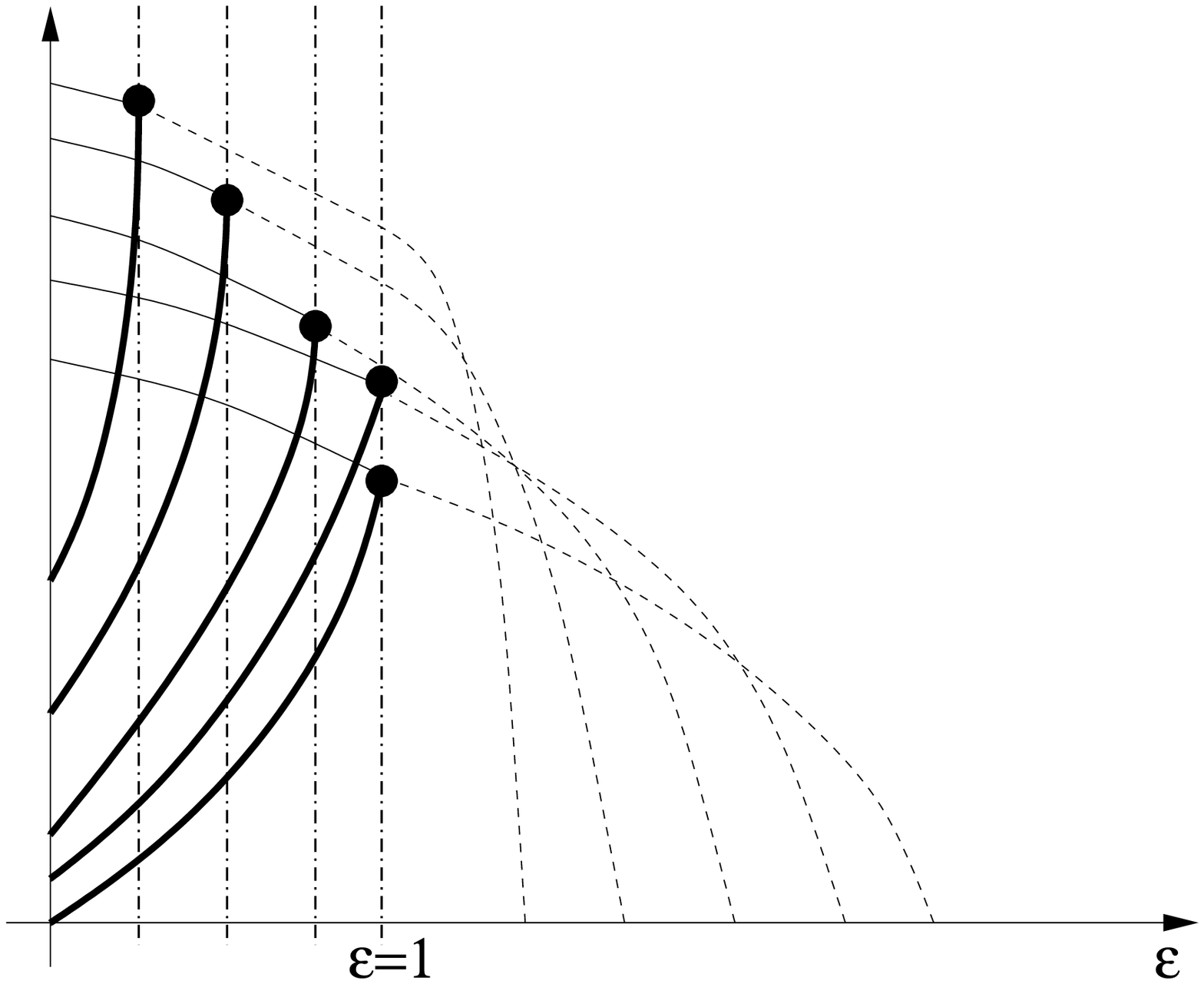}
  \caption{(a) Proposed bifurcation diagram for the one-interface
    solutions in the $N=\infty$ case. We show the primary branch and the
  continuum of saddle-node from the saddle--node bifurcations. In (b)
  we indicate that this structure is then repeated for the
  three-interface, five-interface,... solutions. The solid circles
  represent $\lim_{s \rightarrow 0} \varepsilon_m^{b,s} =
  \varepsilon_m^0  =\varepsilon_m^d$. The bifurcations for the one and
three interface solutions both occur at $\varepsilon=1$.}
\label{fig:general}
\end{center}
\end{figure}

These two conjectures, if true,  would lead to the bifurcation picture
sketched in Figure \ref{fig:general}. In (a) we plot the zero-mean
one-interface branch and have indicated the continuum of
saddle--node bifurcations $\varepsilon_1^s$ that approach the bifurcation at 
$\varepsilon_1^s=\varepsilon_1^d=1$. 
In (b) we indicate the first four branches of the infinite number 
that bifurcate from zero, the branches of associated saddle-node
bifurcations and here we have that $\lim_{s \rightarrow 0}
\varepsilon_m^{b,s} = \varepsilon_m^0  =\varepsilon_m^d$.
In addition our numerical investigation seems to indicate that 
$\varepsilon_1^s=\varepsilon_1^d=1=\varepsilon_2^s=\varepsilon_2^d$.

Finally let us consider the stable solutions - that is the solutions we expect
to see from any simulation. Thus we have for $\epsilon>1$ two stable
solutions, then a region of parameter space with an infinite number of
stable solutions of one and three interface type, then a region of
parameter space with one, three and five interfaces and so on.
In conclusion the diffusion coefficient $\varepsilon$ determines the
number and type of stable solutions.

\bibliography{ref}{}
\bibliographystyle{plain}

\end{document}